\documentclass[a4paper,notitlepage,reqno,twoside,aa]{amsart}
\usepackage{amsmath,amsfonts,amscd,verbatim,amssymb,enumerate,textcomp,mathcomp,stmaryrd,hhline,mathrsfs}
\usepackage{amsthm,graphics}
\usepackage[square]{natbib}
\usepackage[T2A]{fontenc}
\usepackage[cp1251]{inputenc}
\usepackage[russian,english]{babel}

\newcounter{mycounter}
\newcounter{Acounter}

\newtheorem{tth}{Theorem}[section]

\newtheorem{cor}{Corollary}[section]
\newtheorem{lm}{Lemma}[section]

\theoremstyle{definition}
\newtheorem*{defn}{Definition}
\newtheorem*{remn}{Remark}

\theoremstyle{definition}{}
\theoremstyle{definition}{}

\DeclareMathAlphabet{\mathpzc}{OT1}{pzc}{m}{it}

\def \wg{\wedge}

\def \d{\delta}

\def \dk{\delta_k}

\def \G{\Gamma}
\def \Gk{\Gamma_k}
\def \wdw{\wedge\dots\wedge}
\def \e{\varepsilon}
\def \ml{\mathcal}
\def \>{\geqslant}
\def \<{\leqslant}
\def \l{\langle}
\def \r{\rangle}

\def \Lk{\ml{L}_k}
\def \CLk{C_*(\Lk)}

\def \HLk{H_*(\Lk)}

\def \ch{{\rm ch}}
\def \cprime{$'$}
\def \md#1{\equiv -1\,\text{\rm mod}\,#1}
\def \L#1{\mathcal{L}_{#1}}

\def\txt#1{\qquad &\text{\rm #1}\quad}
\def\tx#1{\quad &\text{\rm #1}\quad}

\def\ignore#1{}
\def\rmn#1 {\uppercase\expandafter{\romannumeral#1}}

\def\quat#1{\textquotedblleft #1\textquotedblright}

\title{Note on Laplace operators and homologies\\ of a few Lie subalgebras of $A_1^{(1)}$}

\author{F.~V.~Weinstein}
\address
{Universit\"at Bern, Anatomisches Institut, CH-3012 Bern, B\"uhlstrasse 26, Switzerland.}
\email{Weinstein@ana.unibe.ch}

\begin{document}

\begin{abstract}
In the note we investigate a spectrum of positive self-adjoint operator $\Gk$
(Laplace operator) acting in the external complexes
of some interesting subalgebras $\Lk$ of Lie algebra $A_1^{(1)}$.
We obtain an explicit formula for the action of $\G_k$.
This formula is used to compute the homology of $\Lk$
with trivial coefficients for $k=-1,0,1,2$. In these cases
we show that a spectrum of $\G_k$ is the set of non-negative
integers with a finite multiplicity of each eigenvalue
for $k=-1,0$, infinite one for $k=1,2$. We also found the
generating functions for the multiplicities
of $\G_k$-eigenvalues (in appropriate sense for $\G_1$ and $\G_2$).
\end{abstract}

\maketitle
\markboth{NOTE ON LAPLACE OPERATORS AND HOMOLOGIES}{F.~V.~WEINSTEIN}

The \quat{Laplace operators}, which we consider in the note, appear in the following general situation.

Let $C_*=\bigoplus_qC_q$ be a chain complex over real numbers $\mathbb{R}$
with differential $d$. It is called the \emph{Euclidean complex} if $\dim(C_q)<\infty$,
and on $C_*$ is defined an Euclidean inner product $\l,\r$,
(that is, a non-degenerate positive definite symmetric bilinear form)
so that spaces $C_q$ are pairwise orthogonal.

Let $C^*=\bigoplus_q\mathrm{Hom}_{\mathbb{R}}(C_q,\mathbb{R})=\bigoplus_qC^q$
be the corresponding dual complex. The differential of $C^*$
(that is, a conjugate to $d$ operator) we denote by $\d$.
The inner product defines for each $q$ a canonical isomorphism $C_q\cong C^q$.
Identifying $C_*$ with $C^*$ we may consider operator $\d$ as defined on $C_*$.
Then the linear map
\[
\G:C_*\to C_*,\qquad \G=d\d+\d d
\]
makes sense, and is called \emph{Laplace operator}, or shortly \emph{laplacian}.
This definition evidently extends to a direct sum of the Euclidean complexes.

The efficiency of using the laplacian in homological calculations is based on
three simple claims:
\begin{enumerate}
\item The eigenvalues of $\G$ are non-negative.
\item Let $C^\lambda_*\subset C_*$ be the subspace of $\G$-eigenvectors with eigenvalue $\lambda$.
Then $C_*=\bigoplus_\lambda C^\lambda_*$ is a direct sum of complexes.
\item $H_*(C^\lambda_*)=0$ if $\lambda\neq 0$, and $H^*(C_*)=H_*(C_*)=H_*(C^0_*)=C^0_*$.
\end{enumerate}

\emph{If a chain belongs to the kernel of $\G$ then it is called a harmonic chain of $\G$}.
\begin{proof}
Since $\G$ is a self-adjoint linear operator acting in a finite-dimensional
Euclidian space $C_q$, the first claim follows from a standard fact of linear algebra.
The second one is a corollary of equality $d\G=\G d$.

Let $c\in C^\lambda_q$ and $d(c)=0$. Then $\lambda c=\G(c)=d\d(c)$.
If $\lambda\neq 0$ then $c=d\big(\lambda^{-1}\d(c)\big)$. That is, each cycle of
$C^\lambda_q$ is a boundary. Then $H_q(C^\lambda_*)=0$.
Let $\lambda=0$. Then $0=\l\G(c),c\r=\l\d(c),\d(c)\r+\l d(c),d(c)\r$.
Since our inner product is positive and non-degenerate we conclude that $d(c)=0$,
that is, the differential of complex $C^0_*$ is trivial. Hence $H_q(C^0_*)=C_q^0$.
\end{proof}

If complex $C_*$ has an algebraic origin (e.g., it is an external complex of a Lie algebra)
and its inner product is related with the corresponding algebraic structure,
then not only the harmonic chains, but also the spectral resolution of $\G$ may be a subject of interest.

The laplacian for semisimple finite-dimensional Lie algebras was introduced
in fundamental article by B.~Kostant \cite{MR0142696}.
The subsequent development of Kostant's ideas have impressive applications to Lie algebras
and their homologies (e.g., \cite{MR0387361,MR0414645,MR1923198,MR515625,MR1254731}).

The present note is an elementary investigation of laplacian
for a few naturally defined Lie subalgebras of the simplest
infinite-dimensional Kac-Moody Lie algebra $A_1^{(1)}=\widehat{sl}_2$.
The next two sections contain the main definitions and results summary of the note.

\section{Preliminaries}\label{pref}
\setcounter{mycounter}{0}

Let $sl_2$ be a $3$-dimensional Lie algebra of real $2\times 2$-matrices with zero trace. Take
\[
e_{-1}=
\begin{pmatrix}
0 &0 \\
1 &0
\end{pmatrix},
\quad e_0=\frac{1}{2}
\begin{pmatrix}
1 &0 \\
0 &-1
\end{pmatrix},
\quad e_1=\frac{1}{2}
\begin{pmatrix}
0 &1 \\
0 &0
\end{pmatrix}
\]
as a basis of $sl_2$.
Let $\widehat{sl}_2=\big(sl_2\otimes\mathbb{R}[z^{-1},z]\big)\bigoplus\mathbb{R}p$
be a Lie algebra in which the vectors $p$ and
\[
e_{3h-1}=e_{-1}\otimes z^h,\quad e_{3h}=e_0\otimes z^h,\quad e_{3h+1}=e_1\otimes z^h,\qquad(h=0,\pm 1,\pm 2,\dots),
\]
constitute an $\mathbb{R}$-linear basis of $\widehat{sl}_2$,
and the bracket is expressed by the formulas
\[
[e_a,e_b]=\e_{b-a}e_{a+b},\qquad [p,e_a]=\pi(e_a)e_a,
\]
where $\e_m=-1,0,1$ for $m\md{3}$ respectively, and $\pi(e_a)=(a-\e_a)/3$.
\smallskip

In what follows we assume that $k\>-1$ is an integer.
All linear objects are defined over the field of real numbers $\mathbb{R}$.
The designation $V_*$ means that $V_*=\bigoplus_{q\>0} V_q$
is a graded vector space, where the summands constitute a chain complex.
The homology we interpret as a chain complex with zero differential.

Let \emph{$\Lk$ be a subalgebra of $\widehat{sl}_2$ generated by $e_i$ with $i\>k$},
and $\CLk$ be the \emph{standard complex} of $\Lk$ with coefficients in trivial $\Lk$-module $\mathbb{R}$,
that is, $C_q(\Lk)=\bigwedge^q\Lk$ (see \cite{MR0024908}).

For set $I=\{i_1,\dots,i_q\}$ of the distinct integers $\>k$
a chain of the form $e_I=e_{i_1}\wdw e_{i_q}\in C_q(\Lk)$ is called a $q$-dimensional $k$-\emph{monomial}.
The $k$-monomials with $i_1<\dots<i_q$ constitute a \emph{monomial basis} of $ C_q(\Lk)$.
Differential $d$ of complex $\CLk$ acts by
\[
d(e_{i_1}\wdw e_{i_q})=\sum_{1\< a<b\< q}(-1)^{a+b-1}
\e_{i_b-i_a}e_{i_a+i_b}\wg e_{i_1}\wdw\widehat e_{i_a}\wdw\widehat e_{i_b}\wdw e_{i_q}.
\]
The adjoint actions of $p$ and $e_0$ on $\CLk$ are diagonal on $k$-monomials.
Define the linear functions of \emph{weight} $\omega:\CLk\to\mathbb{Z}$,
and of \emph{degree} $\pi:\CLk\to\mathbb{Z}_{\>0}$ by
\[
e_0(e_I)=(\e_{i_1}+\dots+\e_{i_q})\,e_I=\omega(e_I)\,e_I,\qquad
p(e_I)=\big(\pi(e_{i_1})+\dots+\pi(e_{i_q})\big)\,e_I=\pi(e_I)\,e_I
\]
Let $C_*^{(h)}(\Lk)$ be the subspace of $\CLk$ generated by $k$-monomials $e_I$ such that $\pi(e_I)=h$,
and $C_*^{(w,h)}(\Lk)\subset C_*^{(h)}(\Lk)$ is generated by ones with $\omega(e_I)=w$.
Since $d$ commutes with the actions of $e_0$ and $p$,
\[
\CLk=\bigoplus_{h\>\pi(e_k)}\;C_*^{(h)}(\Lk)=\bigoplus_{(w,h)\in D_k}\;C_*^{(w,h)}(\Lk)
\]
are the direct sums of finite dimensional subcomplexes of $\CLk$;
$D_k$ denotes a set of two-dimensional vectors $(w,h)$,
which may be presented as a sum of the distinct vectors
$\big(\omega(e_a),\pi(e_a)\big)$, $e_a\in\Lk$.
By convention $\dim C_0^{(0,0)}(\Lk)=1$ for any $k$, that is, $(0,0)\in D_k$.
Obviously
\begin{equation}\label{ProdSum}
\prod_{a=k}^\infty\big(1+tu^{\omega(e_a)}x^{\pi(e_a)}\big)=\sum_{(w,h)\in D_k}\sum_{q=0}^\infty\;\dim C_q^{(w,h)}(\Lk)t^qu^wx^h.
\end{equation}

An inner product, defined on the basis of $k$-monomials by formula
\[
\l e_{I_1},e_{I_2}\r=
\begin{cases}
1 &\text{if $I_1=I_2$},\\
0 &\text{if $I_1\neq I_2$},\\
\end{cases}
\]
turns the above decompositions of $\CLk$ into a direct sum of the Euclidian complexes.

Let $\dk$ be an operator on $\CLk$ conjugated to $d$ with respect to the introduced inner product.
Then
\begin{gather*}
\dk(e_i)=\sum\limits_{a+b=i;\ k\< a<b}\;\e_{b-a}\,e_{a}\wg e_{b},\\
\dk(e_{i_1}\wdw e_{i_q})=\sum_{s=1}^q(-1)^{s-1}\;e_{i_1}\wdw\dk(e_{i_s})\wdw e_{i_q}.
\label{E:mm1.2}
\end{gather*}
Denote by $\Gk$ the laplacian of $\CLk$. Then $\text{\rm Ker}(\Gk)=\HLk$.
Obviously complexes $C_*^{(h)}(\Lk)$ and $C_*^{(w,h)}(\Lk)$ are $\G_k$-invariant.
\emph{In what follows we shall identify the space of homologies
$\HLk$ with the space of harmonic chains of $\Gk$}; we also call them the \emph{harmonic chains of $\L{k}$}.

The next terminology and notation will be used: for a linear operator $A$ acting on a linear space $V$,
and complex number $\lambda$, a \emph{$\lambda$-space $V_{\lambda}$ of $A$}
is the eigenspace of $A$ with eigenvalue $\lambda$.
The \emph{$\lambda$-multiplicity of $A$} is the dimension of $V_\lambda$.
\begin{remn}
We have defined the weight $\omega$ as the half of the classically defined one
in  representation theory of $sl_2$. To avoid confusions we included at the end of
note an appendix that contains the basic formulations of finite-dimensional $sl_2$-theory
adopted to our choice of basis in $sl_2$ and definition of weight.
In addition where we introduce the graded $sl_2$-modules, and singular characters of them,
which are our main tools in sections \ref{secLk} and \ref{homol}.
\end{remn}


\section{Summary of results}\label{summary}
Our starting point is an explicit formula for $\Gk$
expressed by means of the adjoint action of $\Lk$.
In the note we mainly consider the cases $k=-1,0,1,2$. Then we show
that the spectrum of $\G_k$ coincides with the set of non-negative integers.
(Remark that the spectrum of $\G_{k>2}$ always includes the irrational numbers.)
Moreover, the multiplicity of each eigenvalue of $\G_k$ is finite for $k=-1,0$ and infinite for $k=1,2$.

The main our goal is an investigation of the spectral resolution of $\G_k$.
Our considerations lead to a simple computing of the (known) homologies of algebras $\Lk$ when $k=-1,0,1,2$,
and to the formulas for multiplicities of eigenvalues of $\G_k$ (in appropriate sense for $k=1,2$).
The material concerning the cases $k=0,1$ should be regarded as well known.

Homologically the most interesting case is $k=2$.
The only familiar to the author computing of $H_*(\L{2})$
is based on a completely different idea (see \cite{MR1254731,Wfb}).
Using the laplacian gives more than just the dimensions of $H_q(\L{2})$.
Namely, since $\L{2}\subset\L{-1}$ is an ideal, and
$\L{-1}/\L{2}\simeq sl_2$, each homology space $H_q(\L{2})$ is an $sl_2$-module.
We show that this module is simple.
Then it is generated by a singular vector, which may be explicitly found.
That gives an elegant formula for the harmonic chains of $\L{2}$.

Section \ref{Resol} concerns with cases $k=0,1$.
Then the eigenvectors of $\G_k$ are $k$-monomials.
Let $C_*^{[w,\lambda]}(\L{k})$ be the subspace of $c\in C_*(\L{k})$
such that $\omega(c)=w$, and $\G_k(c)=\lambda c$.
We use our computing of $H_*(\L{1})$
to prove the Gauss-Jacobi identity. Basing on it we establish the following formulas
\begin{equation*}\label{mintro}
\sum_{\lambda=0}^\infty\dim C_*^{[w,\lambda]}(\L{1})\;x^\lambda=\frac{1}{\Theta(-x,1)}=\prod_{m=1}^\infty\frac{1+x^m}{1-x^m}\;,\quad
\sum_{\lambda=0}^\infty\sum_{w=-\infty}^\infty\dim C_*^{[w,\lambda]}(\L{0})\;u^wx^\lambda=2\,
\frac{\Theta(x,u)}{\Theta(-x,1)},
\end{equation*}
where
\[
\Theta(x,u)=1+2\;\sum_{r=1}^\infty u^rx^{r^2}
\]
is the Jacobi theta-function.
In particular for the $\lambda$-multiplicities of $\G_0$ we get
\[
\sum_{\lambda=0}^\infty\dim C_{*,\lambda}(\L{0})\;x^\lambda=2\,\frac{\Theta(x,1)}{\Theta(-x,1)}
=2\,\prod_{m=1}^\infty\frac{(1+x^{2m-1})^2}{(1-x^{2m-1})^2}\,.
\]

In cases $k=-1,2$, which we consider in Section \ref{homol}, the situation is more complex and interesting. Using
our computing of $H_*\big(\L{2}\big)$, we establish such $\L{2}$-analog of the Gauss-Jacobi identity
\begin{equation*}
\prod_{m=1}^\infty(1-u^{-1}x^m)(1-x^m)(1-ux^m)=
\sum_{w=0}^\infty(-1)^w\;[2w+1]_u\;x^{\frac{w(w+1)}{2}},
\end{equation*}
where $[a]_u=(u^{a/2}-u^{-a/2})/(u^{1/2}-u^{-1/2})$.
One can deduce it from the usual Gauss-Jacobi identity. But our reasonings clearly show that a
natural context of its appearance is the homology of $\L{2}$.
For $u=1$ this identity turns into the famous Jacobi formula for cube of
the Euler function. For $u^3=1,\,u\neq 1$ we obtain
\[
\prod_{m=1}^\infty(1-x^{3m})=\sum_{w=0}^\infty(-1)^w\;\e_{2w+1}\;x^{\frac{w(w+1)}{2}}
\]
that is equivalent to Euler's pentagonal theorem.

Since $\L{2}\subset\L{-1}$ is an ideal and $sl_2\subset\L{-1}$, it follows that $C_*(\L{k})$ for $k=-1,2$
are the $sl_2$-modules.
Each of them may be decomposed into a direct sum of the finite-dimensional $sl_2$-submodules.
Let $S_*^{(w)}(\L{k})\subset C_*(\L{k})$ be the subspace of $sl_2$-singular vectors of dominant weight $w\>0$,
and $S_*^{(w,h)}(\L{k})=S_*^{(w)}(\L{k})\bigcap C_*^{(h)}(\L{k})$.

We shall see that $S_*^{(w,h)}(\L{k})$ is a finite-dimensional $\lambda=\big(h-(-1)^kw(w+1)/2\big)$-space of $\G_k$,
and $C_{*,\lambda}^{(h)}(\L{k})$ is an $sl_2$-submodule generated by $S_*^{(w,h)}(\L{k})$.
Because $h$ is uniquely defined by $w$ and $\lambda$ we may
denote  $S_*^{(w,h)}(\L{k})$ by $S_*^{[w,\lambda]}(\L{k})$. We prove that
\begin{gather*}
\sum_{\lambda=0}^\infty\sum_{w=0}^\infty\dim S_*^{[w,\lambda]}(\L{2})z^wx^\lambda=
\bigg(\frac{1}{1-z}+2\sum_{r=1}^\infty(-1)^r\frac{x^{r(r+1)/2}}{1-zx^r}\bigg)\Theta^{-1}(-x,1),\\
\sum_{\lambda=0}^\infty\sum_{w=0}^\infty\dim S_*^{[w,\lambda]}(\L{-1})z^wx^\lambda
=2\sum_{w=0}^\infty z^w\big(x^{w^2}-x^{(w+1)^2}\big)\;\Theta^{-1}(-x,1).
\end{gather*}
Each side of these identities may be interpreted as \quat{singular characters} of $sl_2$,
that is, as elements of a power series algebra $R(sl_2)[[x]]$
on variable $x$ with coefficients in the representation ring of $sl_2$, where $z^w$ corresponds
to a simple $sl_2$-module of dominant weight $2w+1$ \big(see Appendix\big).
Then the lefthanded sides one can present as infinite products
in $R(sl_2)[[x]]$ of some polynomials \big(see the details in Section \ref{secLk}\big).

As a corollary of the last formula we obtain the generating functions for $\dim C_*^{[w,\lambda]}(\L{-1})$,
and for the $\lambda$-multiplicities of $\L{-1}$. The last one astonishingly coincides
with the one of algebra $\L{0}$:
\[
\sum_{\lambda=0}^\infty\dim C_{*,\lambda}(\L{-1})\;x^\lambda=2\,\frac{\Theta(x,1)}{\Theta(-x,1)}.
\]
(Remark that $\dim C_*^{[w,\lambda]}(\L{2})=\infty$ for all $w,\lambda$.)

In between we also obtain the generating functions for $C_*^{(w,h)}(\Lk)$.
According to \eqref{ProdSum} they may be presented as products
\[
\prod_{a=\frac{k+1}{3}}^\infty(1+u^{-1}x^a)(1+x^a)(1+ux^a)=
\begin{cases}
\sum\limits_{w=0}^\infty(-1)^w[2w+1]_{-u}\;x^{\frac{w(w+1)}{2}}\;\Theta^{-1}(-x,1)&\text{if $k=2$},\vspace{2.5mm}\\
2\sum\limits_{w=0}^\infty [2w+1]_u\,x^{\frac{w(w-1)}{2}}\big(1-x^{2w+1}\big)\;\Theta^{-1}(-x,1)&\text{if $k=-1$}.
\end{cases}
\]

Up to our knowledge the laplacian for $\Lk$ with $k\neq 1$ was not considered earlier.
It is worth to note that $\HLk$ are computed for all $k\>-1$.
For  $k=-1,0$ the result is well known, and easily proved
without using the laplacian. For $k=1$ it is a very particular case of \cite{MR0387361}.
For all $k\>1$ the computation is done in \cite{MR1254731} or \cite{Wfb}.

Our approach to computing $H_*(\L{2})$ may be generalized to computing the homology of
similar to $\L{2}$ subalgebras of the affine Kac-Moody algebras. We hope to consider this problem later.

The author would like to thank D.~Donin for informing on his observation:
the laplacian of $\L{2}$ has an integer spectrum on the space of two-dimensional chains.
A desire to explain it was the initial stimulus for appearing of this note.


\section{Laplacian of algebra $\Lk$}\label{lapk}
\setcounter{mycounter}{0}

A linear endomorphism $T$ of $C_*(\Lk)$ is called a \emph{second order operator} if
\begin{multline*}
T(e_{i_1}\wdw e_{i_q})=\sum_{1\<a<b\<q}(-1)^{a+b-1}\;
T(e_{i_a}\wg e_{i_b})\wg e_{i_1}\wdw\widehat{e}_{i_a}\wdw\widehat{e}_{i_b}\wdw e_{i_q}\\
-(q-2)\sum_{1\<a\<q}e_{i_1}\wdw T(e_{i_a})\wdw e_{i_q}.
\end{multline*}
The next claim is routinely verified.
\stepcounter{mycounter}
\begin{lm}\label{lm_two}
$\G_k:C_*\big(\Lk\big)\to C_*\big(\Lk\big)$ is a second order operator.
\end{lm}
Let $k\>1$, and $e_{-r}$ be an endomorphism of $C_*(\Lk)$,
conjugate to the adjoint action of $e_r$ on $C_*(\Lk)$.
It is a differentiation of the external algebra $\CLk$ of $\Lk$, and
\[
e_{-r}(e_a)=
\begin{cases}
\e_{a+r}e_{a-r}\tx{if} a-r\>k,\\
0\tx{if} a-r<k.
\end{cases}
\]
\stepcounter{mycounter}
\begin{tth}\label{tthlap}
If $k\>1$ then
\begin{equation}\label{eqDk}
\G_k=p+\frac{1}{2}\big(\e_{k+1}^2\,e_0-\sum\limits_{-k<r<k}e_{-r}e_r\big).
\end{equation}
\end{tth}
\begin{proof}
Easily verified that the righthanded side of
equality \eqref{eqDk} is a second order operator on $C_*(\Lk)$.
Therefore in view of Lemma \ref{lm_two} we need only to show that the actions of the both sides of formula \eqref{eqDk}
coincide on one- and two-dimensional $k$-monomials.

The action of $\Gk$ on such monomials may be established
by a direct (but cumbersome for  2-dimensional monomials) calculation.
Namely, let $\lfloor\alpha\rfloor$ be the integer part (floor) of $\alpha\in\mathbb{R}$.
Define
\[
E_k(a)=
\begin{cases}
\big\lfloor\frac{a-2k+2}{3}\big\rfloor\txt{if} a\>2k-2,\\
\; 0\txt{if} a<2k-2.
\end{cases}
\]
Then for $e_a\in\Lk$
\begin{equation*}\label{eq01}
\Gk(e_a)=E_k(a)\;e_a,
\end{equation*}
and for $e_a\wg e_b,\,e_x\wg e_y\in C_2(\Lk),\,(a<b,x<y,a+b=x+y)$
\begin{equation*}\label{eqlmD}
\l\Gk(e_a\wg e_b),e_x\wg e_y\r=
\begin{cases}
E_k(a)+E_k(b)-\e_a\e_b\txt{if}(x,y)=(a,b)\;\;\text{\rm and}\;\;b-a\>k,\\
E_k(a)+E_k(b)+\e_{a-b}^2\txt{if}(x,y)=(a,b)\;\;\text{\rm and}\;\;b-a<k,\\
-\e_{a+x}\e_{b+y}\txt{if either}0<x-a<k\<y-a,\\
\txt{or}0<a-x<k\<b-x,\\
\e_{b-a}\e_{y-x}\txt{if either}a<x\;\;\text{\rm and}\;\;y-a<k,\\
\txt{or}x<a\;\;\text{\rm and}\;\;b-x<k,\\
0\txt{otherwise.}
\end{cases}
\end{equation*}
A verification whether the action of righthanded side of \eqref{eqDk}
on  one- and two-dimensional $k$-monomials
is identical to one of $\G_k$, is a direct calculation as well.
We omit the verification details.
\end{proof}
\stepcounter{mycounter}
\begin{cor}\label{cth1}
\begin{align*}
\G_{-1}&=p+\frac{1}{2}\big(e_{-1}e_1+e_0^2+e_1e_{-1}\big), & \G_0 &=p+\frac{1}{2}\big(e^2_0+e_0\big),\\
\G_1   &=p-\frac{1}{2}\big(e^2_0-e_0\big),                 & \G_2 &=p-\frac{1}{2}\big(e_{-1}e_1+e_0^2+e_1e_{-1}\big).
\end{align*}
\end{cor}
\begin{proof}
The formulas for $\G_1$ and $\G_2$ are the particular cases of equality \eqref{eqDk}.
To prove the formulas for $\G_0$ and $\G_{-1}$ we need only to verify them (as before)
for one- and two-dimensional monomials.

The laplacian $\G_0$ acts on $C_*(\L{0})=C_*(\L{1})\bigoplus e_0\wg C_*(\L{1})$.
Easily checked that for $c\in C_*(\L{1})$
\[
\G_0(c)=\big(e_0^2+\G_1\big)\,c\,,\qquad\G_0(e_0\wg c)=e_0\wg\big(e_0^2+\G_1\big)\,c\,,
\]
what obviously implies the claimed formula for $\G_0$.

The laplacian $\G_{-1}$ acts on $C_*(\L{-1})=C_*(\L{0})\bigoplus e_{-1}\wg C_*(\L{0})$.
Then for $c\in C_*(\L{0})$
\[
\G_{-1}(c)=(e_{-1}e_1+\G_0)\,c\,.
\]
Indeed, as $\d_{-1}(e_a)=\e_{a-1}e_{-1}\wg e_{a+1}+\d_0(e_a)$,
this equality immediately follows for $e_a\in C_1(\L{0})$, as well as for $e_a\wg e_b\in C_2(\L{0})$ since
\begin{multline*}
\G_{-1}(e_a\wg e_b)=d\big(\d_{-1}(e_a)\wg e_b\big)-d\big(e_a\wg\d_{-1}(e_b)\big)+\e_{b-a}\d_{-1}(e_{a+b})=\\
(e_{-1}e_1+\G_{0})\,e_a\wg e_b+
(\e_{b-a}\e_{a+b-1}-\e_{a-1}\e_{b-a-1}-\e_{b-1}\e_{b-a+1})e_{-1}\wg e_{a+b+1},
\end{multline*}
and $\e_{b-a}\e_{a+b-1}-\e_{a-1}\e_{b-a-1}-\e_{b-1}\e_{b-a+1}=0$.
Because $e_0=e_1e_{-1}-e_{-1}e_1$ in universal enveloping algebra of $sl_2$,
the restriction of $\G_{-1}$ to $C_*(\L{0})$
has the form $\G_{-1}=e_{-1}e_1+\G_0=p+(e_{-1}e_1+e_0^2+e_1e_{-1}\big)/2$.
It remains to verify the formula for $\G_{-1}$
for monomials $e_{-1}$ and $e_{-1}\wg e_a$, what is a straightforward testing
\end{proof}


\section{Homology of $\Lk$ and multiplicities of $\G_k$ for $k=0,1$}\label{Resol}
\setcounter{mycounter}{0}

In this section $k=0,1$.
From Corollary \ref{cth1} it follows that a basis of the harmonic chains
of $\G_k$ constitute $k$-monomials $e_I$ such that
\begin{equation*}\label{S01}
\pi(e_I)=\frac{(-1)^{k+1}}{2}\;\omega(e_I)\big(\omega(e_I)+(-1)^k\big).
\end{equation*}

For a fixed $I=\{i_1,\dots,i_q\}$ denote by $q_-,\,q_0,\,q_+$ the quantities of $i'$s,
which are congruent respectively to $-1,0,1\,\text{\rm mod}\,3$. Then
\begin{equation*}\label{eqpm}
\omega(e_I)=q_+ -q_-\,,
\end{equation*}
and for $e_I\in C_q(\Lk)$ evidently
\begin{equation*}
\pi(e_I)\>Q_k(I),\qquad\text{\rm where}\qquad Q_k(I)=\begin{cases}
\frac{q_+(q_+-1)}{2}+\frac{q_-(q_-+1)}{2}+\frac{q_0(q_0-1)}{2} &\mathrm{if}\;\; k=0,\vspace{1.5mm}\\
\frac{q_+(q_+-1)}{2}+\frac{q_-(q_-+1)}{2}+\frac{q_0(q_0+1)}{2} &\mathrm{if}\;\; k=1.
\end{cases}
\end{equation*}

This estimation and the preceding equality show that
\[
Q_0(I)\<-\frac{1}{2}\;(q_+ -q_-)(q_+ -q_-+1),\qquad Q_1(I)\<\frac{1}{2}\;(q_+ -q_-)(q_+ -q_--1).
\]
The inequality for $Q_0(I)$ implies $q_+=q_-$ and thus $Q_0(I)=0$. Then it follows that
$q_0=1$ (since $q>0$). Hence $I=\{0\}$, and $e_0$ is a unique harmonic chain for $\L{0}$.

The inequality for $Q_1(I)$ is equivalent to $\;q_0(q_0+1)\<-2q_+q_-$\,.
Since $q_+,q_-,q_0\>0$ it follows that $q_0=0$, and either $q_+=0$, or $q_-=0$,
which respectively correspond to $I=\{1,4,\dots,3q-2\}$, and $I=\{2,5,\dots,3q-1\}$.
In both cases $\G_1(e_I)=0$. Thus we obtain
\stepcounter{mycounter}
\begin{tth}\label{corL1}
Let $q>0$. Then
\noindent
\begin{enumerate}
\item[\rm{(1)}]
\[
H_q(\L{0})=
\begin{cases}
\mathbb{R}\;e_0 &\text{\rm if\;\;\;$q=1$},\\
0 &\text{\rm otherwise.}
\end{cases}
\]
\item[\rm{(2)}]
The $q$-dimensional harmonic chains $e_1\wg e_4\wdw e_{3q-2}$ and
$e_2\wg e_5\wdw e_{3q-1}$ are the cycles of algebra $\L{1}$.
They represent all $q$-dimensional homological classes
of $\L{1}$. In particular, $\dim H_q(\L{1})=2$ for $q>0$, and
\[
\dim H_q^{(w,h)}(\L{1})=
\begin{cases}
1 &\text{\rm if\;\;\;$w=\pm q,\;\;h=\frac{q^2-q}{2}$},\\
0 &\text{\rm otherwise.}
\end{cases}
\]
\end{enumerate}
\end{tth}

Now let us turn to the multiplicities of $\G_k$.
Consider first algebra $\L{1}$. Corollary \ref{cth1} shows that for $c\in C_*^{(w,h)}(\L{1})$
\begin{equation*}
\G_1(c)=\big(h-w(w-1)/2\big)\,c,
\end{equation*}
that is, $C_*^{(w,h)}(\L{1})$ is an eigenspace of $\G_1$.
Remark that $(w,h)\in D_1$ iff $h-w(w-1)/2\>0$.
Really, the necessity follows because the eigenvalues of $\G_1$ are non-negative.
The sufficiency follows as well since for $\lambda=h-w(w-1)/2\>0$ the vectors
$e_{3\lambda},\,e_1\wg e_4\wdw e_{3(w-2)+1}\wg e_{3(w-1+\lambda)+1}$, and
$e_2\wg e_5\wdw e_{3(w-2)-1}\wg e_{3(w-2+\lambda)-1}$ belong to $C_*^{(w,h)}(\L{1})$
for $w=0, w>0$, and $w<0$ respectively.
Since $h$ is uniquely defined by $\lambda$ and $w$ we obtain
\[
C_*^{[w,\lambda]}(\L{1})=C_*^{(w,\lambda+w(w-1)/2)}(\L{1}).
\]
Therefore equality \eqref{ProdSum} we may rewrite as
\begin{multline}\label{eqt}
(1+tu)\prod_{m=1}^\infty(1+tu^{-1}x^m)(1+tx^m)(1+tux^m)\\
=\sum_{\lambda=0}^\infty\sum_{w=-\infty}^\infty\bigg(\sum_{q=0}^\infty t^q
\dim C_q^{[w,\lambda]}(\L{1})\bigg)\;u^wx^{\lambda+\frac{w(w-1)}{2}}.
\end{multline}
In this formula the sum over $q$ for $t=-1$ is the Euler characteristic of complex $C_*^{[w,\lambda]}(\L{1})$.
By standard claim of homological algebra and Theorem \ref{corL1}(2) it equals to
\[
\sum_{q=0}^\infty\;(-1)^q\dim H_q^{[w,\lambda]}(\L{1})=
\begin{cases}
(-1)^w\qquad &\text{if $\lambda=0$},\\
0\qquad &\text{otherwise.}
\end{cases}
\]
Thus for $t=-1$ equality \eqref{eqt} turns into  relation
\begin{equation}\label{eqGJ}
(1-u)\prod_{m=1}^\infty(1-u^{-1}x^m)(1-x^m)(1-ux^m)=
\sum_{w=-\infty}^\infty (-1)^wu^wx^{\frac{w(w-1)}{2}},
\end{equation}
that is a variant of the \emph{Gauss-Jacobi identity}.
The substitution $u\to u^2x,\,x\to x^2$ turns it into one of its traditional forms
\begin{equation}\label{Jac}
\prod_{m=1}^\infty(1-u^{-2}x^{2m-1})(1-x^{2m})(1-u^2x^{2m-1})=\sum_{w=-\infty}^\infty(-1)^wu^{2w}x^{w^2}.
\end{equation}

\stepcounter{mycounter}
\begin{tth}\label{tthGen1}
For arbitrary $w\in\mathbb{Z}$
\[
\sum_{\lambda=0}^\infty\dim C^{[w,\lambda]}_*(\L{1})\,x^\lambda=\Theta^{-1}(-x,1).
\]
\end{tth}
\begin{proof}
Note first that
\begin{equation}\label{theta}
\prod_{m=1}^\infty\frac{1+x^m}{1-x^m}=\Theta^{-1}(-x,1).
\end{equation}
Really, the substitution $z=\sqrt{-1},\,x\to -x$ in identity \eqref{Jac}
turns the righthanded side of it into $\Theta(-x,1)$,
whereas the lefthanded side turns into
\[
\prod_{m=1}^\infty(1-x^{2m-1})\prod_{m=1}^\infty(1-x^m)=
\prod_{m=1}^\infty\frac{1-x^m}{1-x^{2m}}\prod_{m=1}^\infty(1-x^m)=
\prod_{m=1}^\infty\frac{1-x^m}{1+x^m}.
\]
So \eqref{theta} follows. Let $\Theta^{-1}(-x,1)=\sum_{\lambda=0}^\infty \mu(\lambda)\,x^\lambda$.
Multiplying the result of substitution $u\to -u$ in \eqref{eqGJ} with equality \eqref{theta} we obtain
\begin{equation*}\label{eqt10}
(1+u)\prod_{m=1}^\infty(1+u^{-1}x^m)(1+x^m)(1+ux^m)=
\sum_{\lambda=0}^\infty\sum_{w=-\infty}^\infty\;\mu(\lambda)u^wx^{\lambda+\frac{w(w-1)}{2}}.
\end{equation*}
A comparison of this equality with \eqref{eqt}, where we set $t=1$, completes the proof.
\end{proof}
Now let $k=0$.
Corollary \ref{cth1} shows that for $c\in C_*^{(w,h)}(\L{0})$
\begin{equation*}
\G_0(c)=\big(h+w(w+1)/2\big)\;c.
\end{equation*}
Therefore $C_*^{(w,h)}(\L{0})$ is an eigenspace of $\G_0$.
Obviously $D_0=D_1$, and
\[
C_*^{[w,\lambda]}(\L{0})=\textstyle{C_*^{[w,\lambda-w^2]}(\L{1})\bigoplus e_0\wg C_*^{[w,\lambda-w^2]}(\L{1}).}
\]
Using this decomposition and Theorem \ref{tthGen1} we obtain
\begin{multline*}
\sum_{\lambda=0}^\infty\sum_{w=-\infty}^\infty\dim C_*^{[w,\lambda]}(\L{0})u^w x^\lambda=
2\sum_{w=-\infty}^\infty\;\sum_{\lambda=w^2}^\infty\dim C_*^{[w,\lambda-w^2]}(\L{1})u^w x^\lambda\\
=2\sum_{w=-\infty}^\infty u^wx^{w^2}\sum_{\lambda=0}^\infty\dim C_*^{[w,\lambda]}(\L{1})x^\lambda=
2\;\frac{\Theta(x,u)}{\Theta(-x,1)}\,.
\end{multline*}

\begin{remn}
Similarly to the above proof of formula \eqref{theta} it is easily verified that
\[
\frac{\Theta(x,1)}{\Theta(-x,1)}=\prod_{m=1}^\infty\frac{(1+x^{2m-1})^2}{(1-x^{2m-1})^2}\,.
\]
\end{remn}

\begin{remn}
We obtained the Gauss-Jacobi identity as a corollary of computing the homology of $\L{1}$.
The idea of such a proof (including a proof of much more general Macdonald's identities)
is due to Garland (see \cite{MR0387361,MR0414645}).
There are numerous analytical proofs (see e.g., \cite{MR0011320}).
It is possible to interpret this identity by means of partitions,
and give an elementary bijective proof as it is done in \cite{MR1728812}.
\end{remn}


\section{Standard complex of $\Lk$ for $k\md{3}$}\label{secLk}
\setcounter{mycounter}{0}
This section contains a preparatory material to the calculations in Section \ref{homol}.

Let $k\md{3}$. Then $\Lk\subset\L{-1}$ is an ideal. In particular
on $C_*(\Lk)$ is defined a structure of $sl_2=\{e_{-1},e_0,e_1\}$-module.
\stepcounter{mycounter}
\begin{lm}\label{lm3}
If $k\md{3}$ then $p,\,d,\,\dk $ are the endomorphisms of $sl_2$-module $C_*(\Lk)$.
\end{lm}
\begin{proof}
For $p$ the claim follows because $[p,sl_2]=0$.
Since $\Lk\subset\L{-1}$ is an ideal,
$d$ commutes with the adjoint action of $\L{-1}$ on $C_*(\Lk)$
(this is a general claim).
In particular $d$ commutes with $sl_2$.

The conjugate operator to the adjoint action of $g=a_{-1}e_{-1}+a_0e_0+a_1e_1\in sl_2$
on $\Lk$ is the adjoint action of $\widetilde{g}=a_{-1}e_1+a_0e_0+a_1e_{-1}\in sl_2$.
Then for arbitrary $x,y\in C^*(\Lk)$ we obtain
\[
\l\dk g(x),y\r=\l x,\widetilde{g}d(y)\r=\l x,d\widetilde{g}(y)\r=\l g\dk(x),y\r.
\]
Hence $\dk g=g\dk $ for all $g\in sl_2$, because our inner product is non-degenerate.
\end{proof}

For $sl_2$-module $C^{(h)}_*(\Lk),\,(h>0)$
\emph{let $D(k,h)$ be its set of the dominant weights, $P^{(w,h)}_*(\Lk)$ its isotypic
component of dominant weight $w$, and $S^{(w,h)}_*(\Lk)\subset P^{(w,h)}_*(\Lk)$
be the subspace of singular vectors.}
From Lemma \ref{lm3} and Theorem \ref{sl3}(4) it follows that
\[
C^{(h)}_*(\Lk)=\bigoplus_{w\in D(k,h)}\;P^{(w,h)}_*(\Lk)=
\bigoplus_{w\in D(k,h)}\;\bigoplus_{m=0}^{2w} e_{-1}^m\big(S^{(w,h)}_*(\Lk)\big)
\]
is a direct sum of the isomorphic complexes. Moreover, it follows that
\begin{equation*}\label{DecLk10}
C^{(h)}_{*,\lambda}(\Lk)=\bigoplus_{w\in D(k,h)}\;\bigoplus_{m=0}^{2w} e_{-1}^m\big(S^{(w,h)}_{*,\lambda}(\Lk)\big).
\end{equation*}
In particular, to find the harmonic chains of $\G_k$ it is sufficient to find the spaces of singular harmonic chains
$S^{(w,h)}_{*,0}(\Lk)$ for all $w\in D(k,h)$.

A description of the singular vectors in $C_*(\Lk)$ does not depends on $k$
since $\Lk$ is a direct sum of the adjoint representations of $sl_2$:
\[
\Lk=\bigoplus_{a=(k+1)/3}^\infty M_a,
\qquad\text{where}\qquad M_a=\mathbb{R}e_{3a-1}\oplus\mathbb{R}e_{3a}\oplus\mathbb{R}e_{3a+1}.
\]
Therefore $\CLk$ are isomorphic for different $k\md{3}$ as $sl_2$-modules.

Let $E(a)=\bigoplus_{q=0}^3 M_a^{(q)},\,(M_a^{(q)}=\bigwedge^q M_a)$, be the external algebra of $sl_2$-module $M_a$.
Obviously there is an isomorphism of the $sl_2$-modules
\begin{equation}\label{Decomp}
C_*(\L{k})\cong\bigotimes_{a=(k+1)/3}^\infty E(a)=
{\textstyle{\mathbb{R}\bigoplus}}\bigoplus_{(k+1)/3\<a_1<\dots<a_n}
\Big(\textstyle{M_{a_1}^{(q_1)}\bigotimes\dots\bigotimes M_{a_n}^{(q_n)}}\Big),
\end{equation}
where $q_1,\dots,q_n\in\{1,2,3\}$.
Since each $M_{a}^{(q)}$ is a simple $sl_2$-module of dominant weight $\e_q^2$,
the summands in the last sum are the $sl_2$-modules isomorphic to $V(1)^{\otimes r}$,
where $r=\e_{q_1}^2+\dots+\e_{q_n}^2$.
(As in Appendix, $V(w)$ denotes a simple $sl_2$-module of dominant weight $w\in\mathbb{Z}_{\>0}$.)
Theorem \ref{sl4}(2) supplies an algorithm to find all singular vectors in such modules.

Now we use the singular characters $S$ and notation introduced in Appendix.
It is clear that $S\big(V(1)^{\otimes r}\big)=z^{\otimes r}$.
We can present $z^{\otimes r}$ as a polynomial
\[
z^{\otimes r}=Q_r(z)=z^r+a_{r-1}z^{r-1}+\dots+a_1z+a_0\in R(sl_2).
\]
Consider $E(a)$ as a bigraded $sl_2$-module by setting $\deg(M_{a}^{(q)})=(q,a)$.
Then for the singular characters with values in $R(sl_2)[[t,x]]$ we have
\[
S\big(E(a)\big)=1+ztx^a+zt^2x^{2a}+t^3x^{3a},
\]
\[
S\Big(\textstyle{M_{a_1}^{(q_1)}\bigotimes\dots\bigotimes M_{a_n}^{(q_n)}}\Big)=
t^{q_1+\dots+q_n}x^{a_1+\dots+a_n}Q_{\e_{q_1}^2+\dots+\e_{q_n}^2}(z).
\]

Obviously coefficient $a_w$ of polynomial $Q_r(z)$ equals to
dimension of the singular vectors subspace of isotypic component $\widetilde{V}(w)$
in the primary decomposition of $V(1)^{\otimes r}$.
Let us pass to the singular characters in the both sides of decomposition \eqref{Decomp}.
Then from Lemma \ref{lmA} in $R(sl_2)[[t,x]]$ we obtain the equality
\begin{equation}\label{eqt2}
{\textstyle{\bigotimes\limits_{a=(k+1)/3}^\infty(1+ztx^a+zt^2x^{2a}+t^3x^{3a})}=}
\sum_{(w,h)\in D(k)}\bigg(\sum_{q=0}^\infty\;t^q \dim S_q^{(w,h)}(\L{k})\bigg)\;z^wx^h,
\end{equation}
where $D(k)=\bigcup_{h\>0}D(k,h)$. (By definition $\dim S_*^{(0,0)}(\L{k})=1$,
i.e., $D(k,0)$ is a one-element set.)
\stepcounter{mycounter}
\begin{remn}\label{remRio}
Easily checked that $z^{\otimes r}=Q_r(z)$ is a polynomial, defined by recurrence
\[
Q_r(z)=\frac{z^2+z+1}{z}\,Q_{r-1}(z)-\frac{z+1}{z}\,Q_{r-1}(0),\qquad Q_0(z)=1.
\]
The coefficients of polynomials $Q_r(z)$ are known in combinatorics as the Riordan arrays.
The constant terms of them is a sequence of so called
Motzkin sums (see \cite{MR1327059}, sequence M2587).
\end{remn}

\stepcounter{mycounter}
\begin{remn}\label{rem2.5.8}
For arbitrary $k=3r-1,\,(r\>1)$ it is easy to indicate some
harmonic chains collection of $\Lk$ as follows. The $q$-dimensional chain
$c_q=e_{3r+1}\wg\dots\wg e_{3(r+q-1)+1}\in C_q(\ml{L}_{3r-1})$ is evidently harmonic.
Since $e_0(c_q)=qc_q,\,e_1(c_q)=0$, vector $c_q$
is a singular vector of $sl_2$-module $H_q(\ml{L}_{3r-1})$ of dominant weight $q$.
Then a cyclic $sl_2$-module, spanned by $c_q$ is a $(2q+1)$-dimensional subspace of $H_q(\ml{L}_{3r-1})$.
Because $\G_{3r-1}$ commutes with the action of $sl_2$, the chains
$c_q,e_{-1}(c_q),e_{-1}^2(c_q),\dots,e^{2q}_{-1}(c_q)$ are harmonic.
In the next section we show that they exhaust the harmonic chains of $\L{2}$.
\end{remn}


\section{Homology of $\Lk$ and multiplicities of $\G_k$ for $k=-1,2$}\label{homol}
\setcounter{mycounter}{0}

Assume now $k=-1,2$.
Corollary \ref{cth1} and Theorem \ref{sl3}(4) imply that for $c\in P_*^{(w,h)}(\Lk)$
\[
\G_k(c)=\lambda_k(w,h)\,c,\qquad\text{where}\qquad\lambda_k(w,h)=h-(-1)^kw(w+1)/2.
\]
Therefore $P_*^{(w,h)}(\Lk)$ is a $\lambda_k(w,h)$-eigenspace of $\G_k$.
In particular $s\in S^{(w,h)}_{*,0}(\Lk)$ iff
\begin{equation*}\label{eqlap}
p(s)-\frac{(-1)^k}{2}\;\omega(s)\big(\omega(s)+1\big)s=0
\end{equation*}
(where, of course, $p(s)=hs,\omega(s)=w$). Let $s$ be a linear combination of distinct $k$-monomials $e_I$'s.
Then $\big(\pi(e_I),\omega(e_I)\big)=(w,h)$ for all $e_I$'s.
Since $e_I$'s are linearly independent, each $s=e_I$ satisfies to \eqref{eqlap}.
Thus to find the singular solutions of \eqref{eqlap} first we shall find its monomial solutions,
and then decide, which linear combinations of them are the singular vectors.
So assume that
\[
\pi(e_I)=\frac{(-1)^k}{2}\;\omega(e_I)\big(\omega(e_I)+1\big).
\]
For a fixed $e_I\in C_q(\Lk)$ the next estimation is evident
(notation see at the beginning of Section \ref{Resol}):
\begin{equation*}
\pi(e_I)\>Q_k(I),\qquad\text{\rm where}\qquad Q_k(I)=\begin{cases}
\frac{q_+(q_+-1)}{2}+\frac{q_-(q_--1)}{2}+\frac{q_0(q_0-1)}{2} &\mathrm{if}\;\; k=-1,\vspace{1.5mm}\\
\frac{q_+(q_++1)}{2}+\frac{q_-(q_-+1)}{2}+\frac{q_0(q_0+1)}{2} &\mathrm{if}\;\; k=2.
\end{cases}
\end{equation*}
Applying it to the preceding equation we obtain the inequalities
\[
Q_{-1}(I)\<-\frac{1}{2}\;(q_+-q_-)(q_+-q_-+1),\qquad Q_2(I)\<\frac{1}{2}\;(q_+ -q_-)(q_+ -q_-+1).
\]
The inequality for $Q_{-1}(I)$ gives $q_+=q_-$ and $Q_{-1}(I)=0$.
All such $e_I\in C_*(\L{-1})$ are exhausted by vectors $e_0,\,e_{-1}\wg e_1,\,e_{-1}\wg e_0\wg e_1$, and
only $e_{-1}\wg e_0\wg e_1$ is a singular one.

The inequality for $Q_2(I)$ is equivalent to $q_0(q_0+1)\<-q_-(2q_++1)$.
Then $q_-=q_0=0$, and $\omega(e_I)=q$.
But only $2$-monomial $e_4\wg e_7\wdw e_{3q+1}$ with such properties
satisfies the preceding equation, and it is a singular vector. Thus
\begin{equation}\label{HS}
H_q\big(S^{(w,h)}_*(\L{2})\big)=
\begin{cases}
\mathbb{R}\;e_4\wg e_7\wdw e_{3q+1} &\text{\rm if\;\;\;$w=q,\;\;h=\frac{q^2+q}{2}$,}\\
0 &\text{\rm otherwise.}
\end{cases}
\end{equation}
Collecting all together we obtain
\stepcounter{mycounter}
\begin{tth}\label{corh2}
Let $q>0$. Then
\begin{enumerate}
\item[\rm{(1)}]
\[
H_q(\L{-1})=
\begin{cases}
\mathbb{R}\;e_{-1}\wg e_0\wg e_1 &\text{\rm if\;\;\;$q=3$,}\\
0 &\text{\rm otherwise.}
\end{cases}
\]
\item[\rm{(2)}]
The harmonic chains $e_{-1}^r(e_4\wg e_7\wdw e_{3q+1})\in C_q(\L{2})$,
where $r=0,1,\dots,2q$ are non-zero cycles of algebra $\L{2}$.
They represent all $q$-dimensional homological classes of $\L{2}$. In particular,
$\dim H_q(\L{2})=2q+1$ and
\[
\dim H_q^{(w,h)}(\L{2})=
\begin{cases}
1 &\text{\rm if\;\;\;$-q\<w\<q,\;\;h=\frac{q^2+q}{2}$},\\
0 &\text{\rm otherwise.}
\end{cases}
\]
\end{enumerate}
\end{tth}
Note that $(w,h)\in D(2,h)$ iff $\lambda_2(w,h)=h-w(w+1)/2\>0$.
The necessity follows from non-negativity the eigenvalues of $\G_2$.
The sufficiency follows since for $h-w(w+1)/2\>0$ vector
$e_4\wg e_7\wdw e_{3(w-1)+1}\wg e_{3(w+\lambda)+1}\in\CLk$
is a singular one of weight $w\>0$ and degree $h$.

Therefore for $k=2$ we can rewrite formula \eqref{eqt2} as
\begin{equation}\label{eqt20}
{\textstyle{\bigotimes\limits_{a=1}^\infty(1+ztx^a+zt^2x^{2a}+t^3x^{3a})}=}
\sum_{\lambda=0}^\infty\sum_{w=0}^\infty\bigg(\sum_{q=0}^\infty\;t^q
\dim S_q^{\big(w,\lambda+\frac{w(w+1)}{2}\big)}(\L{2})\bigg)\;z^wx^{\lambda+\frac{w(w+1)}{2}}.
\end{equation}
The sum over $q$ for $t=-1$ in view of \eqref{HS} equals to
\[
\sum_{q=0}^\infty\;(-1)^q\dim H_q\big(S_*^{\big(w,\lambda+\frac{w(w+1)}{2}\big)}(\L{2})\big)=
\begin{cases}
(-1)^w\qquad &\text{if $\lambda=0$},\\
0\qquad &\text{otherwise.}
\end{cases}
\]
Thus for $t=-1$ in $R(sl_2)[[x]]$ we obtain
\begin{equation*}\label{Gauss2}
{\textstyle{\bigotimes\limits_{a=1}^\infty(1-x^a)\big(1-(z-1)x^a+x^{2a}\big)}}=
\sum_{w=0}^\infty\;(-1)^wz^wx^{\frac{w(w+1)}{2}}
\end{equation*}
(a singular analog of the Gauss-Jacobi identity).
Let us apply to the both parts of this identity the map $W_u$ (the definition of $W_u$ see in Appendix).
Then Corollary \ref{corApp} implies the formula
\begin{equation}\label{Gauss2n}
\prod_{m=1}^\infty(1-u^{-1}x^m)(1-x^m)(1-ux^m)=
\sum_{w=0}^\infty(-1)^w\;[2w+1]_u\;x^{\frac{w(w+1)}{2}}.
\end{equation}
which may be considered as an \emph{$\L{2}$-analog of the Gauss-Jacobi identity}.

\stepcounter{mycounter}
\begin{tth}\label{ThMultS2}
Let $S_*^{[w,\lambda]}(\L{2})=S_*^{(w,\lambda+w(w+1)/2)}(\L{2})$. Then
\begin{equation*}
\sum_{w=0}^\infty\sum_{\lambda=0}^\infty\dim S_*^{[w,\lambda]}(\L{2})z^wx^\lambda=
\Bigg(\frac{1}{1-z}+2\sum_{r=1}^\infty(-1)^r\frac{x^{\frac{r(r+1)}{2}}}{1-zx^r}\Bigg)\Theta^{-1}(-x,1).
\end{equation*}
\end{tth}
\begin{proof}
The substitution $u\to -u$ in identity \eqref{Gauss2n} and formula \eqref{theta} imply that
\[
\prod_{a=1}^\infty(1+u^{-1}x^a)(1+x^a)(1+ux^a)
=\Big(1+\sum_{w=1}^\infty(-1)^w[2w+1]_{-u}\;x^{\frac{w(w+1)}{2}}\Big)\Theta^{-1}(-x,1).
\]
By induction easily established that for $w>0$
\[
[2w+1]_{-u}=(-1)^w[2w+1]_u+2\sum_{r=0}^{w-1}(-1)^r[2r+1]_u.
\]
Let us substitute this expression in the previous identity, and then apply
to the both sides of result the map $W_u^{-1}$, which is well defined in this case.
It sends $[2a+1]_u$ to $z^a$, and transforms the lefthanded side into
the lefthanded side of formula \eqref{eqt20} with $t=1$. Therefore we get
\begin{equation}\label{Mult2}
\sum_{w=0}^\infty\sum_{h=0}^\infty\dim S_*^{(w,h)}(\L{2})z^wx^h=
\bigg(1+\sum_{w=1}^\infty\Big(z^w+2(-1)^w\sum_{r=0}^{w-1}(-1)^rz^r\Big)\;x^{\frac{w(w+1)}{2}}\bigg)\Theta^{-1}(-x,1).
\end{equation}
Since $S_*^{(w,h)}(\L{2})=S_*^{[w,h-w(w+1)/2]}(\L{2})$, in order to obtain from the last formula
a generating function for $S_*^{[w,\lambda]}(\L{2})$ we must replace
in it each monomial $z^ax^b$ with monomial $z^ax^{b-a(a+1)/2}$. Thus
\begin{multline*}\label{inter}
\sum_{w=0}^\infty\sum_{\lambda=0}^\infty\dim S_*^{[w,\lambda]}(\L{2})z^wx^\lambda=
\bigg(1+\sum_{w=1}^\infty\Big(z^w+2(-1)^w\sum_{r=0}^{w-1}(-1)^rz^r\;x^{\frac{w(w+1)}{2}-\frac{r(r+1)}{2}}\Big)\bigg)\Theta^{-1}(-x,1)\\
=\sum_{w=0}^\infty z^w\;\Big(1+2\sum_{r=1}^\infty(-1)^rx^{rw+\frac{r(r+1)}{2}}\Big)\Theta^{-1}(-x,1).
\end{multline*}
After summing the appeared geometric progressions we obtain the claimed formula.
\end{proof}
\stepcounter{mycounter}
\begin{tth}\label{ThMultS-1}
Let $S_*^{[w,\lambda]}(\L{-1})=S_*^{(w,\lambda-w(w+1)/2)}(\L{-1})$. Then
\begin{equation*}
\sum_{w=0}^\infty\sum_{\lambda=0}^\infty\dim S_*^{[w,\lambda]}(\L{-1})z^wx^\lambda=
2\sum_{w=0}^\infty z^w\big(x^{w^2}-x^{(w+1)^2}\big)\;\Theta^{-1}(-x,1).
\end{equation*}
\end{tth}
\begin{proof}
From formula \eqref{eqt2} it follows that
\[
\sum_{w=0}^\infty\sum_{h=0}^\infty\dim S_*^{(w,h)}(\L{-1})z^wx^h=
2(1+z){\textstyle{\bigotimes}}\sum_{w=0}^\infty\sum_{h=0}^\infty\dim S_*^{(w,h)}(\L{2})z^wx^h.
\]
Since in $R(sl_2)$
\[
(1+z){\textstyle{\bigotimes}}z^n=
\begin{cases}
\frac{(1+z)^2}{z}\cdot z^n&\text{if $n>0$},\\
1+z&\text{if $n=0$},
\end{cases}
\]
for $f(x,z)\in R(sl_2)[[x]]$ we have
\[
(1+z){\textstyle{\bigotimes}}f(x,z)=\frac{(1+z)^2}{z}\big(f(x,z)-f(x,0)\big)+(1+z)f(x,0)=
\frac{(1+z)^2}{z}f(x,z)-\frac{1+z}{z}f(x,0).
\]
Multiplying the both parts of \eqref{Mult2} in $R(sl_2)$ by $2(1+z)$,
and using the last formula, after some calculation we get
\[
\sum_{w=0}^\infty\sum_{h=0}^\infty\dim S_*^{(w,h)}(\L{-1})z^wx^h=
2\,\sum_{w=0}^\infty z^wx^{\frac{w(w-1)}{2}}(1-x^{2w+1})\;\Theta^{-1}(-x,1).
\]
In order to obtain from this expression a generating function for $S_*^{[w,\lambda]}(\L{-1})$ we must replace
in it each monomial $z^ax^b$ with monomial $z^ax^{b+a(a+1)/2}$.
The result of this substitution is the claimed formula.
\end{proof}


\appendix
\section{Finite dimensional $sl_2$-modules}\label{sl2mod}
\setcounter{mycounter}{0}

We consider a Lie algebra $sl_2$, defined over a zero characteristic field $F$
with basis $\{e_{-1},e_0,e_1\}$, and bracket
\[
[e_0,e_{\pm 1}]=\pm e_{\pm 1},\qquad [e_1,e_{-1}]=e_0.
\]
\begin{defn}
Let $V$ be an $sl_2$-module. We say that $v\in V,\,v\neq 0$ is a
\emph{vector of weight $\lambda\in F$} if $e_0(v)=\lambda v$.
The subspace $V_\lambda=\{v\in V\;|\;e_0(v)=\lambda v\}$
is called a weight subspace of weight $\lambda$.
\end{defn}

\begin{defn}
A weight vector $v\in V$ is called a \emph{singular vector}, if $e_1(v)=0$.
A \emph{dominant weight} of $V$ is a weight of any singular vector from $V$.
\end{defn}
The structure of finite-dimensional $sl_2$-modules
is described by the next two theorems (see \cite{MR499562}).
\stepcounter{mycounter}
\begin{tth}\label{sl1}
\noindent
\begin{enumerate}
\item[\rm{(1)}] Any finite-dimensional $sl_2$-module is isomorphic to a direct sum of simple $sl_2$-modules.
\item[\rm{(2)}] Every finite-dimensional $sl_2$-module has a unique singular vector up to collinearity.
\item[\rm{(3)}] Every simple finite-dimensional $sl_2$-module is generated by a singular vector.
\item[\rm{(4)}] The simple finite-dimensional $sl_2$-modules, generated by the singular vectors of same weight, are isomorphic.
\end{enumerate}
\end{tth}
\stepcounter{mycounter}
\begin{tth}\label{sl3}
Let $V(w)$ be a finite-dimensional $sl_2$-module with a singular vector
$v$ of weight $w$. Then
\begin{enumerate}
\item[\rm{(1)}] $2w\in\mathbb{Z}_{\>0}$.
\item[\rm{(2)}] The $sl_2$-module $V(w)$ is simple.
\item[\rm{(3)}] Let $v_p=e^p_{-1}(v)$. Then $\{v_0,v_1,\dots,v_{2w}\}$ is a basis of $V(w)$
and $v_p=0$ for $p>2w$.
\item[\rm{(4)}] The operator $e_{-1}e_1+e_0^2+e_1e_{-1}\in U(sl_2)$ acts on $V(w)$ as multiplication by $w(w+1)$.
\end{enumerate}
\end{tth}
We need also the following
\stepcounter{mycounter}
\begin{tth}{\rm (Clebsch-Gordan)}\label{sl4}
Let $w_1,w_2\in\frac{1}{2}\mathbb{Z}_{\>0}$.
Then
\noindent
\begin{enumerate}
\item[\rm{(1)}] There is an isomorphism of the $sl_2$-modules
\[
\textstyle{V(w_1)\bigotimes V(w_2)\cong V(w_1+w_2)\bigoplus
V(w_1+w_2-1)\bigoplus\dots\bigoplus V\big(|w_1-w_2|)}.
\]
\item[\rm{(2)}] Let $v_1\in V(w_1),v_2\in V(w_2)$ be the singular vectors, and $0\<p\<2w_2$.
Then
\[
\sum_{i=0}^p(-1)^i\frac{(2w_2-p+i)!(2w_1-i)!}{(2w_2-p)!(2w_1)!p!(p-i)!}\;\textstyle{e^i_{-1}(v_1)\bigotimes e^{p-i}_{-1}(v_2)}
\]
is a singular vector in $V(w_1)\otimes V(w_2)$ of weight $w_1+w_2-p$.
\end{enumerate}
\end{tth}

\begin{defn}
A $\mathbb{Z}$-module, freely generated by symbols $z^w$,
$w\in\frac{1}{2}\mathbb{Z}_{\>0}$ with multiplication
\[
\textstyle{z^{w_1}\bigotimes z^{w_2}=z^{w_1+w_2}+z^{w_1+w_2-1}+\dots+z^{|w_1-w_2|}}
\]
is called the \emph{representation ring of $sl_2$} and denoted by $R(sl_2)$.
\end{defn}
For a finite-dimensional $sl_2$-module $V$ let $\widetilde{V}(w)\subset V$
be a submodule, generated by the singular vectors of weight $w$.
Then by Theorem \ref{sl1} there is a unique decomposition
\[
V=\bigoplus_{w\in\mathbb{N}}\widetilde{V}(w)
\]
that is called a \emph{primary decomposition of $V$},
the summands of which $\widetilde{V}(w)$ are called the \emph{isotypic components of $V$}.
Define $S(V)\in R(sl_2)$ by
\[
S(V)=\sum_{w\in\mathbb{N}}m(V,w)z^w,\qquad\text{where}\qquad
m(V,w)=
\begin{cases}
\dim\widetilde{V}(w)/\dim V(w)&\text{if $\widetilde{V}(w)\neq\{0\}$},\\
0&\text{otherwise}.
\end{cases}
\]
Since primary decomposition is unique, $S$ sends isomorphic $sl_2$-modules to the same element of $R(sl_2)$.
Together with Theorem \ref{thChar} this shows that a finite-dimensional $sl_2$-module $V$
up to isomorphism is uniquely defined by $S(V)$.
\begin{defn}
Let $\mathbb{N}^r,\;(r\in\mathbb{Z}_{\>0})$ be a semigroup of vectors
$\alpha=(\alpha_1,\dots,\alpha_r)$ with the non-negative integer coordinates.
An \emph{$r$-graded $sl_2$-module} is a direct sum of the finite-dimensional $sl_2$-modules
$V=\bigoplus_{\alpha\in\mathbb{N}^r}V_\alpha$.
If $v\in V_\alpha,\,v\neq 0$ then we say that $v$ has degree $\alpha$.
A \emph{tensor product of the $r$-graded $sl_2$-modules $V_1,V_2$} is an $r$-graded $sl_2$-module
\[
{\textstyle{V_1\bigotimes V_2}}=\bigoplus_{\alpha\in\mathbb{N}^r}{\textstyle{(V_1\bigotimes V_2)_\alpha}},
\qquad\text{where}\qquad{\textstyle{(V_1\bigotimes V_2)_\alpha}}
=\bigoplus_{\alpha_1+\alpha_2=\alpha}\textstyle{V_{\alpha_1}\bigotimes V_{\alpha_2}}.
\]
\end{defn}
\begin{defn}
Let $R(sl_2)[[x]]=R(sl_2)[[x_1,\dots,x_r]]$ be a ring of the power series
on variables $x_1,\dots,x_r$ with coefficients in $R(sl_2)$.
A \emph{singular character of $r$-graded $sl_2$-module $V$ with value in $R(sl_2)[[x]]$} is a power series
\[
S(V)=\sum_{\alpha\in\mathbb{N}^r}S(V_\alpha)\;x^\alpha\in R(sl_2)[[x]],\qquad (x^\alpha=x_1^{\alpha_1}\dots x_r^{\alpha_r}).
\]
\end{defn}
\stepcounter{mycounter}
\begin{lm}\label{lmA}
Let $V_1,V_2$ be the $r$-graded $sl_2$-modules. Then
\[
\textstyle{S(V_1\bigotimes V_2)=S(V_1)\bigotimes S(V_2)},
\]
where symbol $\bigotimes$ at the righthanded side means the product in $R(sl_2)[[x]]$.
\end{lm}
\begin{defn}
Let $V=\bigoplus_\lambda V_\lambda$ be the
decomposition of a finite-dimensional $sl_2$-module $V$ into a direct sum of the weight subspaces.
(By theorems \ref{sl1} and \ref{sl3} such decomposition exists and is unique.)
The \emph{Weil character of $V$} is the Laurent polynomial
\[
\ch_V(u)=\sum_\lambda \dim(V_\lambda)u^\lambda\in\mathbb{Z}[u^{-1/2},u^{1/2}].
\]
\end{defn}
The Weil characters are the same for the isomorphic $sl_2$-modules.
From Theorem \ref{sl3} easily follows that
\[
\ch_{V(w)}(u)=[2w+1]_u,\qquad\text{where}\qquad
[a]_u=\frac{u^{a/2}-u^{-a/2}}{u^{1/2}-u^{-1/2}},\quad(a\in\mathbb{Z}_{\>0}).
\]
\stepcounter{mycounter}
\begin{tth}\label{thChar}
Let $W_u:R(sl_2)\to\mathbb{Z}[u^{-1/2},u^{1/2}]$ be a $\mathbb{Z}$-linear map such that
$W_u(z^w)=[2w+1]_u$. Then $W_u$ is an isomorphism
of ring $R(sl_2)$ onto subring of $\mathbb{Z}[u^{-1/2},u^{1/2}]$, linearly generated by $[a]_u$,
where $a$ runs over the set of non-negative integers.
\end{tth}
This theorem and Lemma \ref{lmA} imply the following assertion.

\stepcounter{mycounter}
\begin{cor}\label{corApp}
Define a ring homomorphism
\[
W_u:R(sl_2)[[x]]\to \mathbb{Z}[u^{-1/2},u^{1/2}][[x]]\qquad\text{\rm by}\qquad
W_u\big(P(z)x^\alpha\big)=W_u\big(P(z)\big)x^\alpha,
\]
where $P(z)\in R(sl_2)$. Then for $r$-graded $sl_2$-modules $V_1,V_2$,
\[
\textstyle{W_u\big(S(V_1\bigotimes V_2)\big)=W_u\big(S(V_1)\big)\cdot W_u\big(S(V_2)\big)}.
\]
\end{cor}
\begin{remn}
It should be noted that Corollary \ref{corApp}
is easily extended to the infinite tensor products of the $r$-graded $sl_2$-modules.
\end{remn}


\end{document}